\documentclass{IEEEtran4PSCC}

\usepackage[cmex10]{amsmath}
\usepackage[font=small]{subfig} 
\usepackage{stfloats}
\usepackage{url}
\usepackage[normalem]{ulem}
\usepackage{graphicx}
\usepackage{amsfonts,amssymb}
\usepackage{cite}
\usepackage[english]{babel}
\usepackage[utf8]{inputenc}
\usepackage{color}
\usepackage{booktabs}
\usepackage{algorithm}
\usepackage{algpseudocode}
\usepackage{tikz}
\usetikzlibrary{arrows,shapes,snakes}

\DeclareMathOperator{\kmeans}{\mathnormal{k}-means}
\DeclareRobustCommand{\bs}[1]{\ensuremath{\boldsymbol{#1}}}
\DeclareRobustCommand{\bm}[1]{\ensuremath{\boldsymbol{#1}}}

\tikzstyle{grid}=[thick,draw=black,fill=black!3!white]
\tikzstyle{tick}=[thick]
\tikzstyle{trajec}=[very thick]

\definecolor{colora}{rgb}{0.75,0.25,0.25}
\definecolor{colorb}{rgb}{0.65,0.75,0.25}
\definecolor{colorc}{rgb}{0.25,0.75,0.45}
\definecolor{colord}{rgb}{0.25,0.45,0.75}
\definecolor{colore}{rgb}{0.65,0.25,0.75}

\title{Continuous Time Multi-stage Stochastic Reserve and Unit Commitment}

\begin{document}

\author{\IEEEauthorblockN{K\'{a}ri Hreinsson, Bita Analui and Anna Scaglione}\IEEEauthorblockA{ECEE, Arizona State University \\ \{kari.hreinsson,bita.analui,anna.scaglione\}@asu.edu }}

\maketitle

\begin{abstract} 
In this paper we introduce a continuous time multi stage stochastic optimization for scheduling generating units, their commitment, reserve capacities and their continuous time generation profiles in the day-ahead wholesale electricity market. Our formulation approximates the solution of a  variational problem, in which the balance, generation capacity and ramping constraints are in continuous time. Due to the greater accuracy of our representation of ramping events this approach improves the system reliability and lowers the real-time cost. 
\end{abstract}

\thanksto{This work was funded in part by the Advanced Research Projects
Agency-Energy (ARPA-E), U.S. Department of Energy, under Award Number
DWS1103. The views and opinions of authors expressed herein do not
necessarily state or reflect those of the United States Government or any
agency thereof.}

\vspace*{-1em}
\section{Introduction and Motivation}

There is a vast literature dealing with stochastic versions of the unit commitment (UC) and of the security constrained unit commitment (SCUC) problems. 
Currently, uncertainty in power systems is managed by scheduling reserve capacity in advance to compensate for errors in (net-)load forecasts. 

Stochastic Unit Commitment (SUC) appears to be the appropriate tool to capture the exogenous uncertainties directly into the decision process. 
In fact, SUC ensures that a feasible solution exists for all considered scenarios, and that the expected cost is minimized.
The most commonly studied versions of SUC formulations are two-stage SUC problems, pioneered by \cite{Wiebking1977}. 
The \textit{multi-stage} stochastic unit commitment (MSUC) is an extension of the two-stage SUC where decisions are taken sequentially at certain time instants. 
While the MSUC formulation allows for smoother boundary conditions on the commitment variables, it has received less attention than the two-stage SUC because the complexity is often prohibitive. 
First in \cite{Takriti&etal1996} and then in a series of follow up work (see e.g.~\cite{Carpentier&etal1996}, \cite{Nowak&Rom2000}, \cite{Birge&Shiina2004}, \cite{Papa&Oren&ONil2011} and the references therein) many authors worked on curbing the MSUC computational complexity. 
It is natural to use the more accurate representation of the uncertainties in SUC to optimally oversee a more economic commitment of reserves.
The authors of \cite{Bouffard&etal2008} have proposed a scenario based approach for a short-term scheduling problem in the presence of wind power generation. More specifically, their proposed approach focuses on the lack of accurate prediction for wind power generation and the allocation of reserve capacities from Hydro-Thermal generation units to secure the system against immense variability of wind power generation. In a similar stream, our work  provides a novel solution to overcome the existing challenges in committing units and scheduling reserves, but our focus is also on improving the accuracy beyond the conventional stochastic optimization frameworks cited above.

We observe that, while representing the decision uncertainty more accurately leads to greater reliability, one of the aspects that the works cited above ignore is the inadequacy in representing the continuous time scenarios and 
all of the literature cited above fail in capturing the continuous time nature of some actions, such as ramping needs and inter-hourly variations. 
For example, in locations where solar power is significant, early morning and evening hours have pronounced curvatures, due to the mismatch between demand peaks and the solar production. 
  
The main contribution of this paper is to show a viable numerical approximation of the stochastic variational problem that power system operators try to address in their day-ahead commitment and reserve decisions. We leverage our prior work in \cite{Parvania&Scaglione2016} focused on the deterministic continuous time unit commitment.
More specifically, we introduce the continuous time multi-stage reserve and unit commitment (CT-MSRUC) problem, in which an underlying load scenario tree is used to decide the baseline day-ahead dispatch, commitment and reserve capacity, considering continuous time generation trajectories as part of the decision variables.
Our approach allows to better schedule the reserve capacity and power, not only because the formulation accounts for the future expected real-time cost in dispatching them, but also because it
provides a more accurate representation of the possible future inter-hourly ramping needs. 

The paper is organized as follows. 
In Section \ref{sec:scenario-tree} we introduce our approximation of continuous-time trajectories and scenario trees of net-load.
The formulation of the CT-MSRUC follows in Section \ref{sec:MS-CoSRUC}.
The formulation is further shown to reduce to a framework similar to the conventional discrete time representation of the generation trajectory, by relaxing specific constraints associated with the time-continuity and using a first order polynomial.
This relaxation (which we refer to as the discrete-time MSRUC --- DT-MSRUC) is used as a benchmark to assess the performance of the proposed framework in Section \ref{sec:numerical}, where we show that the CT-MSRUC leads to greater system reliability and lowers operating cost.

\section{General Polynomial Interpolation of Electric Net-Load Process on a Scenario Tree}\label{sec:scenario-tree}
To approximate continuous net-load trajectories, we utilize the Bernstein polynomial basis:
\begin{equation}
  b_{k,n}(x) = \binom nk x^k (1-x)^{n-k} \Pi(x), \quad k\in[0,n]
\end{equation}
where $\Pi(x)$ is the rectangular function and $\Pi(x) = 1$ for $0\leq x \leq 1$.
We define the vector of polynomials of degree $n$ as $\mathbf{b}_{n}(x)$. 
The Bernstein polynomial operator $\mathfrak{B}_n$ takes a function $f(x) \in [0,1]$ and maps it into an $n$th order polynomial as $\mathfrak{B}_n f(x) = \textstyle \sum_{k=0}^nf_k b_{k,n}(x)$,
an approximation of $f(x)$ (if Lipschitz continuous) with maximum error $O(1/n)$ when $f_k=f(k/n)$. The coefficients $f_k$ are called {\it control points}.

The uncertainty considered in this paper comes from net-load, which we consider a continuous time, continuous space stochastic process $\Xi(t)$. To make the problem computationally tractable, the $\sigma$-algebra of this process must be approximated by a discrete and finite state space stochastic process, whose probability model is truncated over the finite horizon of the optimization. 
A scenario tree $\mathcal{T}=\{\mathcal{V},\mathcal{E}, \mathbb{P};\bm{\xi}\}$ represents the basic structure that is used in multi-stage stochastic optimization to represent the gradual unfolding of information through filtration $\mathfrak{F}$ of the underlying stochastic process. 
The tree is a directed graph from the root node $v_0$, (node/vertex in graph theory sense) of the tree to its leaf nodes, through edges $(v_-,v)\in \mathcal{E}$, where $v_-$ is the parent of node $v$, and $v, v_- \in \mathcal{V}$.
The scenario tree is constructed to meet some minimum approximation error, relative to a metric that is typically a surrogate for minimizing the optimality loss due to the inaccuracy of the simplified decision model. 
In the conventional construction of a tree (cf.~\cite{Pflug&Pichler2015}), each node of the tree has an associated discrete time net-load sample, and the representation of the statistics is nodal, that is: each node $v\in \mathcal{V}$ of the tree has a certain probability $\pi_v$ that is the marginal probability of that particular sample for the quantized process at that stage. 
A better interpretation of the approximate scenario tree for a continuous time process is that realizations are mapped into piece-wise constant sample paths.

We extend this notion by applying the Bernstein operator of order $n$ to each interval mapping the sample net-load paths into a process whose realizations are random piece-wise polynomial functions of order $n\geq 0$. 
For the resulting process, we can construct the corresponding filtration $\mathfrak{F}$ by representing the sub-$\sigma$ algebras of the set of random coefficients. 
More specifically, for all nodes $v \in \mathcal{V}$, we define the set of its direct successor nodes by $v_+ \in \mathcal{C}(v)$ and assume that $t_{v_-}\leq t < t_{v}$ represents the continuous time between any two consecutive nodes $v_{-}~\text{and}~v$.
We use $\pi_v=\pi_{v_{-},v}$ to represent the joint probability of the free parameters in the polynomial trajectory and use the notation $\bm{\xi}\sim \mathcal{T}$ to refer to the complete construction.
We divide the time between consecutive nodes $[t_{v_-}, t_v]$ into $n+1$ arbitrary subintervals and define a vector of \textit{control points}:
\begin{equation}
  {\bm \xi}_{v_-,v}= \left[\xi^{(0)}_{v_-,v},\ldots,\xi^{(n-1)}_{v_-,v}, \xi^{(n)}_{v_-,v}\right]^{\mathrm{T}} \label{proc}
\end{equation}
to represent the coefficients of net-load in the function space of polynomials of degree $n$.
The continuous net-load scenario approximation in each edge's time interval is:
\begin{equation}
\xi_{v_-,v} (t)=\mathbf{b}_{n}\left(\frac{t-t_{v_-}}{t_v-t_{v_-}}\right) {\bm \xi}_{v_-,v},~~~
t_{v_-} \leq t < t_{v} 
\end{equation}
where $\left(\frac{t-t_{v_-}}{t_v-t_{v_-}}\right)$ scales the continuous time index to fit into the time interval corresponding to the edge $(v_-,v)$. 
This approximation
has several useful properties, one of which is that $\xi_{v_-,v}(t_{v_-}) =\xi_{v_-,v}^{(0)}$ and $\xi_{v_-,v}(t_v) =\xi_{v_-,v}^{(n)}$. 
In order to maintain $C^0$ continuity across edges of the load tree, it is sufficient to enforce that the control points match at the endpoints. 
Another useful property of Bernstein polynomials is that the coefficients of an $n^\text{th}$ degree polynomial's derivative are the finite differences of the original coefficients: $\dot{b}_{k,n}(t) = n(b_{k-1,n-1}(t)-b_{k,n-1}(t))$.
Therefore, a linear matrix ${\bf M}_{n+1 \times n}$, which is rectangular-bidiagonal, can be introduced to translate the derivatives of ${\bf b}_n(t)$ into the same family of polynomials of degree $n-1$: 
\begin{align}
  \dot {\bf b}_n(t) = n{\bf M}{\bf b}_{n-1}(t),
\end{align}
Consequently, the elements of $\dot{\bm \xi}_{v_-,v}$ can be expressed as the linear combination of the control points of ${\bm \xi}_{v_-,v}$ as ${\dot{\bm \xi}_{v_-,v} = n{\bf M} {\bm \xi}_{v_-,v}}$.
Specifically, at the end points: 
\begin{align}
  \dot{\bm \xi}_{v_-,v} = \left[n\left(\xi^{(1)}_{v_-,v} - \xi^{(0)}_{v_-,v}\right),\ldots, n\left(\xi^{(n)}_{v_-,v} - \xi^{(n-1)}_{v_-,v}\right)\right]^{\mathrm{T}} \label{der_proc}
\end{align}
Note that neither vector \eqref{proc} nor \eqref{der_proc} are nodal variables, but carry information across the edge $(v_{-},v)$. Figure \ref{fig:genramp} shows the relationship between a continuous (in this case generation) function and it's coefficients, along with it's derivative.

We construct the scenario tree by recusively applying the well known $k$-means clustering algorithm to a bundle of net-load sample paths.
We consider $H$ stages (hours), and have a $\mathbb{R}^{H\times(n+1)}$ dimensional sample space of polynomial (spline) coefficients. The sample space is partitioned into Voronoi cells that are constrained to satisfy the structure imposed by filtration $\mathfrak{F}$ and the piece-wise continuity of polynomial functions. 
Considering $L$ net-load sample paths, we approximate each with a piece-wise continuous polynomial of degree $n$ using the Bernstein basis. 
The vector of coefficients for each path and each interval are stored in the tensor $\bs{\Xi}^s \in \mathbb{R}^{L\times H \times (n+1)}$ spanning the $H$ stages.
The total number of tree nodes is fixed ahead of time as $|\mathcal{V}|=N+1$, and the number of nodes per stage determined by the entries of the non decreasing $1\times H$ vector $\bs{c}$. 
The edges ${\cal E}$ are allocated optimally while quantizing $\bs{\Xi}^s$.  
The structure of the optimum tree is stored in an ancestry (parent) vector $\bs{a}$ such that node $v$ is a child of node $a^{1}_{v}$, and grandchild of node $a^{2}_{v}$.
The non-anticipativity property is enforced by finding an appropriate mapping between each scenario $\bs{\xi}^s_\ell \in \bs{\Xi}^s$ and the tree that ties the past values to the future samples. 
To summarize, for the scenario tree construction, we map our selection of sample discrete-time net-load trajectories to their continuous time approximations, and then merge paths until we get our desired tree structure.
Figure \ref{fig:trees} shows trees constructed for our numerical comparison.

\section{The Continuous Time Multi-Stage Reserve and Unit Commitment Problem (CT-MSRUC)}\label{sec:MS-CoSRUC}
The major differences between a conventional UC formulation and our CT-MSRUC are (a) we consider continuous load and generation trajectories, approximated through Bernstein coefficients, and (b) the stochastic UC considers multiple future scenarios.
The stochastic nature of the formulation means only a handful of decision variables are defined w.r.t.~time (i.e.~have subscript $t$), instead, most are defined in \emph{nodal} form, and thus use subscript $v$ to identify that they refer to tree node $v$ or edge $(v_-,v)$.

\subsection{Considerations Regarding Optimizing in Continuous-Time}
Ideally, in the continuous-time UC formulation, generation, ramping, commitment, start-up and shut-down decision variables $\mathrm{x}^g(t)$, $\dot{\mathrm{x}}^g(t)$, $\mathrm{y}^g(t)$, $\overline{\mathrm{s}}^g(t)$ and $\underline{\mathrm{s}}^g(t)$, for all generation units $g\in \mathcal{G}$, may vary continually throughout time $t$. 
This provides the system with ultimate flexibility to balance the load in an optimal way. 
Now that the continuous load process and its derivative on the tree are identified, we will extend the idea into the unit commitment formulation for the generation trajectories. 
We introduce:
\begin{align}
{\bf x}^g_{v_-,v} &= [x^{g(0)}_{v_-,v}, x^{g(1)}_{v_-,v}, x^{g(2)}_{v_-,v}\ldots,x^{g(n-1)}_{v_-,v}, x^{g(n)}_{v_-,v}]^{\mathrm{T}} \label{G_coef}\\
\dot{{\bf x}}^g_{v_-,v}& =[\dot{x}^{g(0)}_{v_-,v}, \dot{x}^{g(1)}_{v_-,v},\dot{x}^{g(2)}_{v_-,v}\ldots,\dot{x}^{g(n-1)}_{v_-,v}]^{\mathrm{T}} \label{R_coef}
\end{align}
where ${\bf x}^g_{v_-,v}$ and $\dot{{\bf x}}^g_{v_-,v}$ are the continuous-time polynomial coefficients for generation and ramping trajectories respectively. 
It is clear that the elements of vector $\dot{{\bf x}}^g_{v_-,v}$ can be expressed as a linear combination of elements of generation, $\dot{\bf{x}}^g_{v_-,v} = n{\bf M}{\bf{x}}^g_{v_-,v}$, and the corresponding continuous-time generation and ramp functions are: 
\begin{align}
  \mathrm{x}^g_{v_-,v}(t) &= {\bf b}_n\bigg(\frac{t-t_{v_-}}{t_{v}-t_{v_-}}\bigg){\bf x}^g_{v_-,v} \label{generation}\\
  \dot{\mathrm x}^g_{v_-,v}(t) & = n{\bf M}{\bf b}_{n-1}\bigg(\frac{t-t_{v_-}}{t_{v}-t_{v_-}}\bigg)\dot{\bf x}^g_{v_-,v}\label{ramping}
\end{align}
We assume that the commitment and therefore start-up and shut-down variables are constant within each interval $t_{v_-}\leq t < t_v$ and equal to the commitment, start-up and shut-down decisions at the \emph{end} of the interval. 
The reason being that the control point at the end of the interval $[t_{v_-}, t_v]$, which we conventionally denote by $y^g_v$, $\overline{s}^g_v$ and $\underline{s}^g_v$ respectively, carries all the information of the edge $(v_-,v)$.

Another useful property of Bernstein polynomials states that the entire generation and ramping trajectories for an edge $(v_-,v)$ are contained within the convex hull of their control points ${\bf x}^g_{v_-,v}$ and $\dot{{\bf x}}^g_{v_-,v}$ respectively. 
Thus, we can bound generation dispatch and ramping for $g\in\mathcal{G}, v\in\mathcal{V}^+ = \mathcal{V}\smallsetminus\{0\}$ (all tree nodes excluding the root), and $t\in[t_{v_-}, t_v]$:
\begin{align}
  \min\{{\bf{x}}^g_{v_-,v}\} &\leq {\mathrm x}^g_{v_-,v}(t) \leq \max\{{\bf{x}}^g_{v_-,v}\} \\
  \min\{{\bf{\dot{x}}}^g_{v_-,v}\} &\leq \dot{\mathrm x}^g_{v_-,v}(t) \leq \max\{{\bf{\dot{x}}}^g_{v_-,v}\}
\end{align}%

\subsection{Design of the combined reserve and power market}
Conventional unit commitment problems are solved several hours ahead of the time horizon they apply to, and should provide generator owners with production schedules that include expected dispatch, commitment as well as any reserve power obligations. 
Solving such a problem while considering several possible load outcomes (mapped to a tree) naturally yields several distinct dispatch solutions that may or may not include variations in unit commitment, depending on how the problem is constrained. 
The outcome is a set of decisions, some with respect to time, but others with respect to the different paths of the tree, which is hard to interpret as a market decision. 
This is why we propose to merge the commitment of reserves together with the traditional procurement of generation to meet the base-load, framing all of these decisions as part of our stochastic optimization variables. 

Beside the continuous generation trajectories, the following summarizes the difference between the current operational norm and our formulation, as seen from the perspective of a generator owner:\\
1) The continuous generation and commitment \emph{schedule} is the solution that corresponds to the most likely path on the scenario tree. 
    This is denoted by the spline $x^{g,\text{schedule}}(t)$ and the vector $y^{g,\text{schedule}} \in \{0,1\}^H$ for generation and commitment respectively, where $H$ is the number of hours.\\
2) The variable $\overline{y}^g \in \{0,1\}^H$ indicates at which hour $h$ unit $g$ \emph{may} be committed.
      As such, $\overline{y}_h^g = 1$ indicates that for one or more of the paths of the scenario tree that cross hour $h$, the unit is committed, and the generator owner needs to anticipate this.
      Note that a commitment profile will always honor minimum on/off constraints.\\
3) The continuous reserve capacity of units is denoted by $\underline{r}^{g}(t)$ and $\overline{r}^{g}(t)$, generator $g$ should anticipate to be operating in the range $[x^{g,\text{schedule}}(t) - \underline{r}^{g}(t), x^{g,\text{schedule}}(t) + \overline{r}^{g}(t)]$ at time $t$.\\
4) We split payments for commitment into (a) payment for the possibility of being committed during a particular hour $\overline{Y}^g$ and (b) the real-time payment for actually being committed $Y^g$.

Note that as we commit reserves it would be natural to also include N-1 contingencies, but for simplicity we omit this in the CT-MSRUC. The extension is straightforward.

\subsection{Program Formulation}
To reduce complexity, the formulation does not consider transmission constraints, assuming a single-bus system. 
We assume the problem spans $H$ stages (hours), defining the set $\mathcal{H} \in \{1,2,\ldots,H\}$.
We define the set of all scenario tree nodes $\mathcal{V}$, the set of non-root nodes $\mathcal{V}^+ = \mathcal{V}\smallsetminus\{0\}$, the set of nodes at stage $h$ as $\mathcal{Z}(h)$, and the preceding stage ${h-1=h_-}$. 
We assume the root node to be the only member of $\mathcal{Z}(0)$.
We use the subscript $h$ to refer to the time instant at the end of stage $h$, which pertains to decisions made in the interval $t_{h_-} \leq t < t_h$.
The edge between consecutive nodes $v_-$ and $v$ spans the time $t_{v_-} \leq t < t_v$.
As in the conventional UC formulations, we will consider this interval equal to 1 hour.
To make our notation more compact, we consider the case of cubic Bernstein polynomials ($n=3$).
The generalization to an arbitrary $n$ is straightforward. 
The following sections write out the nodal unit commitment formulation.

\subsubsection{Continuity Constraints}
We enforce $C^0$ and $C^1$ continuity across the edge polynomials. For all $g\in\mathcal{G}$:
\begin{align}
	x_{v_-,v}^{g(n)} = x_{v,v_+}^{g(0)} \quad \forall v\in\mathcal{V}^+, v_{+} \in \mathcal{C}(v), v \not\in \mathcal{Z}(H) \label{eq:continuity1}\\
    \dot{x}_{v_-,v}^{g(n-1)} = \dot{x}_{v,v_+}^{g(0)} \quad \forall v\in\mathcal{V}^+, v_{+} \in \mathcal{C}(v), v \not\in \mathcal{Z}(H) \label{eq:continuity2}
\end{align}

\subsubsection{Reserve Capacity Constraints}
We define the continuous reserve capacity in terms of edge variables as well as the aforementioned hourly reserve capacity indicators.
The edge-wise continuous reserve coefficients $\hat{\mathrm{r}}_{v_-,v}^{g(i)}$ and $\check{\mathrm{r}}_{v_-,v}^{g(i)}$ define the necessary reserve capacity required of the units committed on that edge, while $\overline{\mathrm{r}}_{h_-,h}^{g(i)}$ and $\underline{\mathrm{r}}_{h_-,h}^{g(i)}$ define the continuous hourly reserve requirements considering all edges $v \in \mathcal{Z}(h)$ spanning $t_{v_-} \leq t < t_v$.
For all $g\in\mathcal{G}$:
\par\nobreak\vspace*{-1em}
{\small
\begin{align}
	{\bf \overline {r}}_{h_-,h}^{g} \geq {\bf x }_{v_-,v}^{g} + {\bf{\hat{r}}}_{v_-,v}^{g} - {\bf x}^{g,\text{schedule}}_{h_-,h}  \quad \forall h\in\mathcal{H}, v\in\mathcal{Z}(h) \label{resv1}\\
  {\bf \underline{r}}_{h_-,h}^{g} \geq {\bf x}^{g,\text{schedule}}_{h_-,h}- {\bf{x}}_{v_-,v}^{g} + {\bf{\check{r}}}_{v_-,v}^{g}\quad \forall h\in\mathcal{H}, v\in\mathcal{Z}(h) \label{resv2}
\end{align}
}%
Here, we assume the continuous functions $\underline{r}^g(t), \overline{r}^g(t)$ and $x^{g,\text{schedule}}(t)$ are modelled using splines, in the same fashion as the edge-wise generation $x^g(t)$, thus (note that $t_h - t_{h_-} = 1$):
\begin{align}
  \underline{r}^g(t)       &= {\bf b}_n(t-t_{h_-}) {\bf \underline{r}}_{h_-,h}^g \quad t_{h_-} \leq t \leq t_h\\
  \overline{r}^g(t)        &= {\bf b}_n(t-t_{h_-}) {\bf \overline{r}}_{h_-,h}^g \quad t_{h_-} \leq t \leq t_h\\
  x^{g,\text{schedule}}(t) &= {\bf b}_n(t-t_{h_-}) {\bf x}_{h_-,h}^{g,\text{schedule}} \quad t_{h_-} \leq t \leq t_h
\end{align}
The binary $y_v^g$ denotes whether a unit is committed along the preceding edge $(v_-,v)$ and it's relationship with the hourly commitment indicator $\overline{y}_h^g$ is:
\begin{align}
	\overline{y}_h^g \geq y_v^g \quad\forall h\in\mathcal{H},v\in\mathcal{Z}(h),g\in\mathcal{G}
\end{align}
\subsubsection{Balance Constraints}
We maintain a balance of generation and load by ensuring the coefficients are equal:
\begin{align}
	\textstyle \sum_{g\in \mathcal{G}} {\bf x}^g_{v_-,v}-{\bm \xi}_{v_-,v}=0 \label{C_balance}\quad \forall v\in \mathcal{V}^+
\end{align}
Furthermore, we denote the estimated root mean square error between $\bs{\xi}_{v_-,v}$ and the bundle of scenarios it approximates as $\bs{\epsilon}_{v_-,v}$. 
We can require reserves to cover a certain range around the centroid, proportional to the mean square error:
\begin{align}
	\!\!\!\! \textstyle \sum_{g\in \mathcal{G}} ({\bf{x}}^g_{v_-,v} + {\bf{\hat{r}}}_{v_-,v}^g) - {\bm{\xi}}_{v_-,v} - \rho\, {\bm{\epsilon}}_{v_-,v} &=0 \quad\forall v\in\mathcal{V}^+\\
	\!\!\!\!\textstyle \sum_{g\in \mathcal{G}} ({\bf{x}}^g_{v_-,v} - {\bf{\check{r}}}_{v_-,v}^g) - {\bm{\xi}}_{v_-,v} + \rho\, {\bm{\epsilon}}_{v_-,v} &=0 \quad\forall v\in\mathcal{V}^+ 
\end{align}
For an unbiased centroid $\bs{\xi}$, $\bs{\epsilon}$ corresponds to the standard deviation of the sample paths. Essentially $\rho$ allows us to put a margin equal to a multiple of the estimated conditional standard deviation, which can be used as a bound for the probability of deviating from the $\bs{\xi}$.

\subsubsection{Generation and Ramping Limits}
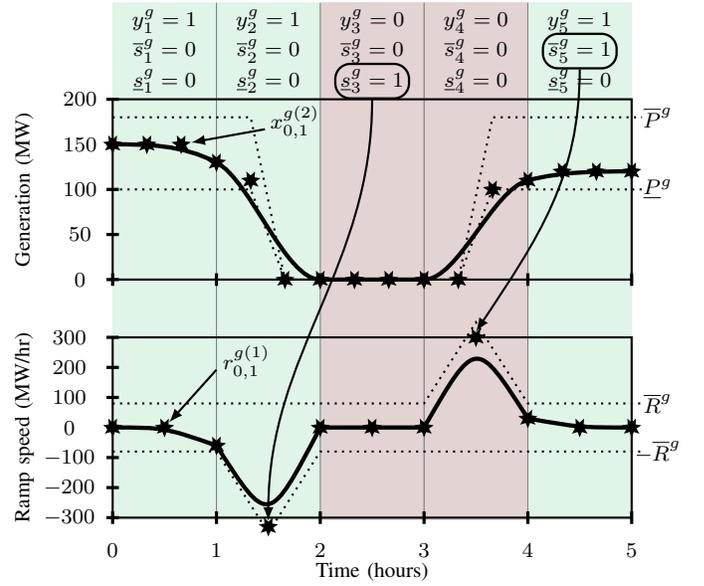
\begin{figure}
  \begin{tikzpicture}[baseline,xscale=0.92,yscale=0.8]
  \footnotesize
  
  \tikzstyle{coeff} = [fill=black,star,star points=7,star point ratio=0.5,minimum size=0.2];

  \draw[draw=none,fill=colorc!80!white,opacity=0.2] (0,-4) rectangle (1.5,4.6);
  \draw[draw=none,fill=colorc!80!white,opacity=0.2] (1.5,-4) rectangle (3,4.6);
  \draw[draw=none,fill=colora!60!black,opacity=0.2] (3,-4) rectangle (4.5,4.6);
  \draw[draw=none,fill=colora!60!black,opacity=0.2] (4.5,-4) rectangle (6,4.6);
  \draw[draw=none,fill=colorc!80!white,opacity=0.2] (6,-4) rectangle (7.5,4.6);
  \draw[draw=black!50!white] (1.5,-4) -- (1.5,4.6);
  \draw[draw=black!50!white] (3.0,-4) -- (3.0,4.6);
  \draw[draw=black!50!white] (4.5,-4) -- (4.5,4.6);
  \draw[draw=black!50!white] (6.0,-4) -- (6.0,4.6);

  \begin{scope}[xscale=1.5,yscale=0.015]

    \draw[thick] (0,0) -- (5,0) -- (5,200) -- (0,200) -- cycle;

    \node[rotate=90] (ylabel) at (-0.85,100) {Generation (MW)};
    \foreach \x in {0,1,2,3,4,5}
      { \draw[very thick] (\x,7) -- (\x,-7); }
    \foreach \y in {0, 50, 100, 150, 200}
      { \draw[very thick] (.06,\y) -- (-.06,\y) node[xshift=-10] {$\y$}; }

    \draw[ultra thick] (0,150) .. controls (0.33,150) and (0.66,150) ..  (1,130); 
    \draw[ultra thick] (1,130) .. controls (1.33,110) and (1.66,0)   ..  (2,0); 
    \draw[ultra thick] (2,0)   .. controls (2.33,0)   and (2.66,0)   ..  (3,0); 
    \draw[ultra thick] (3,0)   .. controls (3.33,0)   and (3.66,100) ..  (4,110); 
    \draw[ultra thick] (4,110) .. controls (4.33,120) and (4.66,120) ..  (5,120); 
    \draw[thick, dotted] (0,180) -- (1.33,180) -- (1.66,0) -- (3.33,0) -- (3.66,180) -- (5.1,180) node[xshift=6] {$\overline{P}^g$};
    \draw[thick, dotted] (0,100) -- (1.33,100) -- (1.66,0) -- (3.33,0) -- (3.66,100) -- (5.1,100) node[xshift=6] {$\underline{P}^g$};

    \node[coeff] at (0.00,150) {}; \node[coeff] at (0.33,150) {}; \node[coeff] (pgex) at (0.66,150) {}; \node[coeff] at (1.00,130) {}; \node[coeff] at (1.33,110) {}; \node[coeff] at (1.66,0) {}; \node[coeff] at (2.00,0) {}; \node[coeff] at (2.33,0) {}; \node[coeff] at (2.66,0) {}; \node[coeff] at (3.00,0) {}; \node[coeff] at (3.33,0) {}; \node[coeff] at (3.66,100) {}; \node[coeff] at (4.00,110) {}; \node[coeff] at (4.33,120) {}; \node[coeff] at (4.66,120) {}; \node[coeff] at (5.00,120) {};
  \end{scope}

  \begin{scope}[yshift=-70,xscale=1.5,yscale=0.005]
    \node[rotate=90] (ylabel) at (-0.85,0) {Ramp speed (MW/hr)};
    \node[rotate=0]  (xlabel) at (2.5,-480) {Time (hours)};
    \foreach \x in {0,1,2,3,4,5}
      { \draw[very thick] (\x,-279) -- (\x,-321) node[yshift=-10] {$\x$}; }
    \foreach \y in {-300, -200, -100, 0, 100, 200, 300}
      { \draw[very thick] (.06,\y) -- (-.06,\y) node[xshift=-14] {$\y$}; }

    \draw[thick] (0,-300) -- (5,-300) -- (5,300) -- (0,300) -- cycle;
    \draw[ultra thick] (0,0)   .. controls (0.5,0)    ..  (1,-60); 
    \draw[ultra thick] (1,-60) .. controls (1.5,-330) ..  (2,0); 
    \draw[ultra thick] (2,0)   .. controls (2.5,0)    ..  (3,0); 
    \draw[ultra thick] (3,0)   .. controls (3.5,300)  ..  (4,30); 
    \draw[ultra thick] (4,30)  .. controls (4.5,0)    ..  (5,0); 
    \node[coeff] at (0.0,0) {}; \node[coeff] (rgex) at (0.5,0) {}; \node[coeff] at (1.0,-60) {}; \node[coeff] (lr) at (1.5,-330) {}; \node[coeff] at (2.0,0) {}; \node[coeff] at (2.5,0) {}; \node[coeff] at (3.0,0) {}; \node[coeff] (ur) at (3.5,300) {}; \node[coeff] at (4.0,30) {}; \node[coeff] at (4.5,0) {}; \node[coeff] at (5.0,0) {};
    \draw[thick, dotted] (0,80)  -- (3.0,80) -- (3.5,350) -- (4.0,80) -- (5.1,80) node[xshift=6] {$\overline{R}^g$};
    \draw[thick, dotted] (0,-80) -- (1.0,-80) -- (1.5,-350) -- (2.0,-80) -- (5.1,-80) node[xshift=6] {$-\overline{R}^g$};
  \end{scope}
  \node at (0.75,4.3) {$y_1^g = 1$}; 
  \node at (0.75,3.8) {$\overline{s}_1^g = 0$}; 
  \node at (0.75,3.3) {$\underline{s}_1^g = 0$};

  \node at (2.25,4.3) {$y_2^g = 1$}; 
  \node at (2.25,3.8) {$\overline{s}_2^g = 0$}; 
  \node at (2.25,3.3) {$\underline{s}_2^g = 0$};

  \node at (3.75,4.3) {$y_3^g = 0$}; 
  \node at (3.75,3.8) {$\overline{s}_3^g = 0$}; 
  \node[draw=black,inner sep=2pt,rounded corners=2mm,line width=0.3mm] (us) at (3.75,3.3) {$\underline{s}_3^g = 1$};

  \node at (5.25,4.3) {$y_4^g = 0$}; 
  \node at (5.25,3.8) {$\overline{s}_4^g = 0$}; 
  \node at (5.25,3.3) {$\underline{s}_4^g = 0$};

  \node at (6.75,4.3) {$y_5^g = 1$}; 
  \node[draw=black,inner sep=2pt,rounded corners=2mm,line width=0.3mm] (os) at (6.75,3.8) {$\overline{s}_5^g = 1$}; 
  \node at (6.75,3.3) {$\underline{s}_5^g = 0$};

  \draw[thick,-latex] (us.south)[out=270,in=90] to (lr.north);
  \draw[thick, -latex] (os.south)[out=270,in=60] to (ur.north);

  \draw[thick, latex-] (pgex) -- (2.2,2.5) node[xshift=0.4cm,yshift=0.1cm] {$x_{0,1}^{g(2)}$};
  \draw[thick, latex-] (rgex) -- (1.5,-1.5) node[xshift=0.4cm,yshift=0.1cm] {$r_{0,1}^{g(1)}$};

\end{tikzpicture}
  \caption{The Bernstein coefficients of cubic spline interpolation for generation and ramping (starred), along with the commitment status and turn on/off indicators. Note how a turn-on or turn-off indicator relaxes the ramping constraint (dotted line) for the center ramping coefficient of the preceding hour.}\label{fig:genramp}
  \vspace*{-1.5em}
\end{figure}%
First off, we define the relationship between generation and ramping:
\begin{align}	
	\dot{{\bf{x}}}^g_{v_-,v} = n{\bf{M}}{\bf{x}}^g_{v_-,v} \quad \forall v\in\mathcal{V}^+, g\in\mathcal{G}
\end{align}
Now, for a unit that is offline for a certain edge (hour) $v_-$ and online for the consecutive edge $v$, we can not simply force all the coefficients ${\bf{x}}_{v_-,v}$ to be zero and ${\bf{x}}_{v,v_+}$ to be non-zero, as this violates the continuity constraints. 
Instead, we allow units to turn-on and turn-off as late as possible during the preceding hour. 
This translates into the following constraints for all units $g\in\mathcal{G}$, assuming $n=3$:
\par\nobreak\vspace*{-1em}
{\small
\begin{align}
	x_{v_-,v}^{g(i)}      + \hat{r}  _{v_-,v}     ^{g(i)} &\leq \overline {P}^g y_v^g \quad \forall v\in\mathcal{V}^+, i\in\{0,1\}\label{eq:pmax1}\\
	x_{v_-,v}^{g(i)}      - \check{r}_{v_-,v}     ^{g(i)} &\leq \underline{P}^g y_v^g \quad \forall v\in\mathcal{V}^+, i\in\{0,1\}\\
  x_{a^2_v,a^1_v}^{g(i)} + \hat{r}  _{a^2_v,a^1_v}^{g(i)} &\leq \overline {P}^g y_v^g \quad \forall v\in\mathcal{V}, v \not\in\mathcal{Z}(1), i\in\{2,3\}\\
  x_{a^2_v,a^1_v}^{g(i)} - \check{r}_{a^2_v,a^1_v}^{g(i)} &\leq \underline{P}^g y_v^g \quad \forall v\in\mathcal{V}, v \not\in\mathcal{Z}(1), i\in\{2,3\}\\
  x_{v_{-},v}^{g(i)}    + \hat{r}  _{v_{-},v}   ^{g(i)} &\leq \overline {P}^g y_v^g \quad \forall v\in\mathcal{Z}(H), i\in\{2,3\}\\
  x_{v_{-},v}^{g(i)}    - \check{r}_{v_{-},v}   ^{g(i)} &\leq \underline{P}^g y_v^g \quad \forall v\in\mathcal{Z}(H), i\in\{2,3\}
\end{align}}%
recalling the ancestry vector $\bs{a}$ introduced in Section \ref{sec:scenario-tree} and the convention $a^0_v := v$.
Generalizing this to any $n$ involves replacing the sets $\{0,1\}$ with $\{0,\ldots,n-C\}$ and $\{2,3\}$ with $\{n-C+1,\ldots,n\}$ where $C$ is the level of continuity ($0$ for $C^0$ and so on).
The ramping constraints become inter-hourly, but we must ensure that they are relaxed for on-off transitions. For all $g\in\mathcal{G}$:
\par\nobreak\vspace*{-1em}
{\small
\begin{align}
  -\overline{R}^g &\leq \dot{x}_{v_{-},v}^{g(i)} \leq \overline{R}^g \quad \forall v\in\mathcal{V}^+, i\in\{0,2\}\label{eq:ramp1}\\
 -\overline{R}^g\! \! -\!  n\overline{P}^g \underline{s}_{v}^g &\leq \dot{x}_{a^2_v,a^1_v}^{g(1)} \leq \overline{R}^g\! \!  +\!  n\overline{P}^g\overline{s}_{v}^g \quad \forall h\in\{2,\ldots,H\}, v\in\mathcal{Z}(h) \label{eq:ramp-relaxed}\\
  \overline{R}^g &\leq \dot{x}_{v_{-},v}^{g(1)} \leq \overline{R}^g \quad\forall v\in\mathcal{Z}(H) \label{eq:ramp2}
\end{align}
}%
Again here, to generalize this to any $n$ involves applying \eqref{eq:ramp-relaxed} and \eqref{eq:ramp2} to coefficient $i = n-C$, and leaving the rest ($\{0,\ldots,n-1\} \smallsetminus \{n-C\}$) to \eqref{eq:ramp1}.
Figure \ref{fig:genramp} attempts to visually explain constraints \eqref{eq:pmax1}-\eqref{eq:ramp2}.

\subsubsection{Minimum On and Off Constrains}
We formulate conventional minimum on-off constraints for all $g\in\mathcal{G}$, where $O^g_n$ and $O^g_f$ denote the minimum on and off time respectively:
\begin{align}
  \overline{s}_v^g - \underline{s}_v^g &= y_v^g - y_{v_-}^g \quad \forall v\in\mathcal{V}^+\\
  y_v^g     &\geq \textstyle\sum_{u=0}^{\min(O_n^g,h)-1} \overline{s}^g_{a^u_v}  \quad\forall h\in\mathcal{H}, v\in\mathcal{Z}(h)\\
  1 - y_v^g &\geq \textstyle\sum_{u=0}^{\min(O_f^g,h)-1} \underline{s}^g_{a^u_v} \quad\forall h\in\mathcal{H}, v\in\mathcal{Z}(h)
\end{align}

\subsubsection{Objective Function}
Integrating Bernstein polynomials is straight-forward (again assuming $t_v-t_{v_-}=1$):
\begin{equation*}
  \int_{t_{v_-}}^{t_v} C^g {\mathrm x}(t) \, dt = C^g \int_{t_{v_-}}^{t_v} {\bf{b}}_n(t-t_{v_-}) {\bf{x}}_{v_-,v}^g = \frac{C^g}{n+1} \bs{1} \cdot {\bf{x}}_{v_-,v}^g
\end{equation*}
In this case, the objective of the stochastic unit commitment is minimizing the expected cost of generation and the reserve capacity payment.
The objective becomes:
\begin{equation}
  \begin{split}
    \min \sum_{h\in\mathcal{H}} \sum_{g\in\mathcal{G}} \left(\frac{1}{n+1}\sum_{i=0}^n \left[\overline{R}^g\overline{r}_{h_-,h}^{g(i)} + \underline{R}^g\underline{r}_{h_-,h}^{g(i)}\right] + \overline{Y}^g \overline{y}_h^g\right)\\
      + \sum_{v\in\mathcal{V}^+} \pi_v \sum_{g\in\mathcal{G}} \left(Y^gy_v^g + \overline{S}\overline{s}^g_v + \underline{S}\underline{s}^g_v + \frac{X^g}{n+1}\sum_{i=0}^n x_{v_-,v}^{g(i)}\right) 
  \end{split} \label{eq:objective}
\end{equation}
where $\pi_v$ is the probability of node $v$, $\overline{R}$ and $\underline{R}$ are up- and down-spinning reserve costs, $\overline{S}$ and $\underline{S}$ are startup and shutdown costs, $X$ is the linear generation cost term, $Y$ is the constant (commitment) cost term and $\overline{Y}$ is the payment for possibly being committed for a particular hour, a reserve payment of sorts. The time-based segment of \eqref{eq:objective} denotes up-front payments for reserves and anticipated operation range, while the nodal based segment is the expected real-time cost. Without information about real-time prices, we simply substitute real-time prices with the day-ahead bids, as the best indicator of what will happen in real-time. However, we wish to remark that having a higher real-time cost will only make the numerical comparison more favorable to our framework.

\section{Numerical Results}\label{sec:numerical}

\subsection{Load Scenario Tree}
\begin{figure}
  \centering
  \begin{tikzpicture}[xscale=1,yscale=1]
  \footnotesize
  \def\x{2.6}
  \def\y{1.0}
  
  \tikzstyle{datapoint} = [fill=black,diamond,inner sep=0pt,minimum size=5pt];
  \draw[grid] (0,0) rectangle (3*\x,3*\y);
  \begin{scope}[xscale=\x,yscale=\y]
    \input{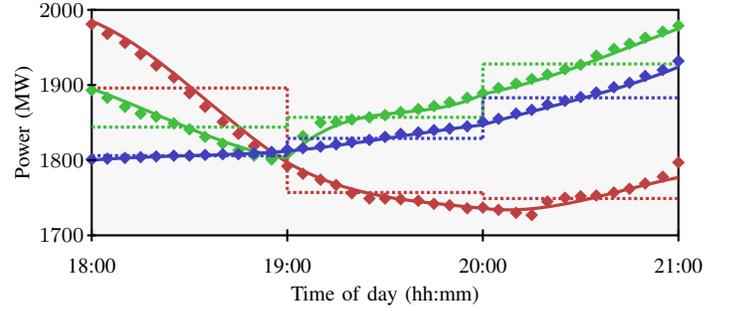}
  \end{scope}

  \foreach \yy in {0,1,2,3} { \draw[tick] (-0.06,\yy*\y) -- (+0.06,\yy*\y); }
  \node at (-0.4, 0*\y) {$1700$};
  \node at (-0.4, 1*\y) {$1800$};
  \node at (-0.4, 2*\y) {$1900$};
  \node at (-0.4, 3*\y) {$2000$};
  \foreach \xx in {0,1,2,3} { \draw[tick] (\xx*\x,-0.08) -- (\xx*\x,+0.08); }
  \node at (0*\x, -0.4) {18:00};
  \node at (1*\x, -0.4) {19:00};
  \node at (2*\x, -0.4) {20:00};
  \node at (3*\x, -0.4) {21:00};

  \node[rotate=90] at (-0.9,1.5*\y) {Power (MW)};
  \node[] at (1.5*\x,-0.8) {Time of day (hh:mm)};

\end{tikzpicture}
  \vspace*{-0.8em}
  \caption{The spline approximation of three training scenarios for two hours. The sample data is denoted by diamonds, the zero order piecewise constant spline is dotted while the $C^1$ continuous cubic spline is solid. We see that the cubic spline is already an extremely good approximation of the sample path.}\label{fig:splines}
  \vspace*{-0.5em}
\end{figure}
\begin{figure}
  \subfloat[The zero-order spline scenario tree.]{\includegraphics[width=0.95\columnwidth]{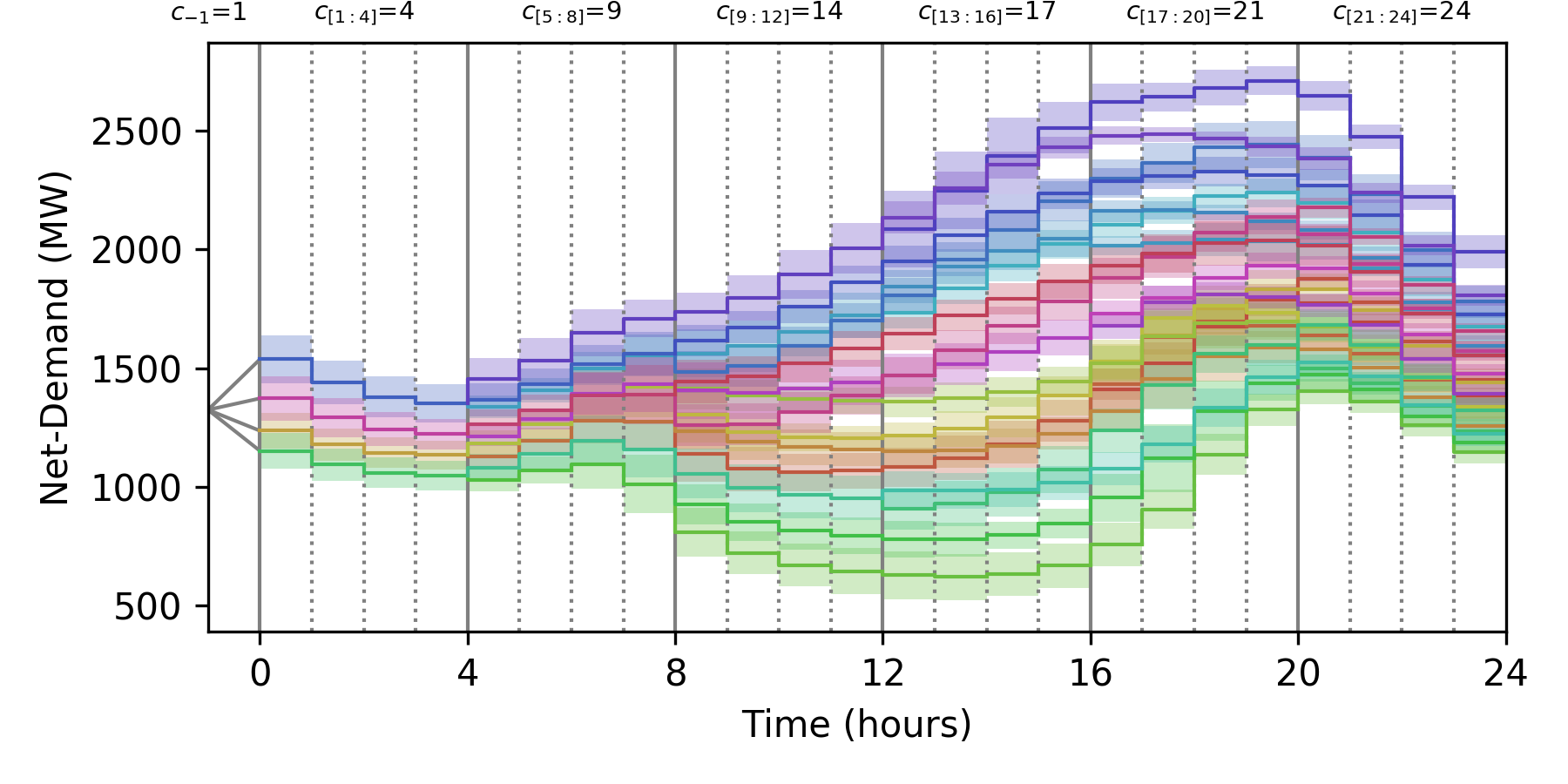}}\\
  \subfloat[The cubic spline scenario tree.]{\includegraphics[width=0.95\columnwidth]{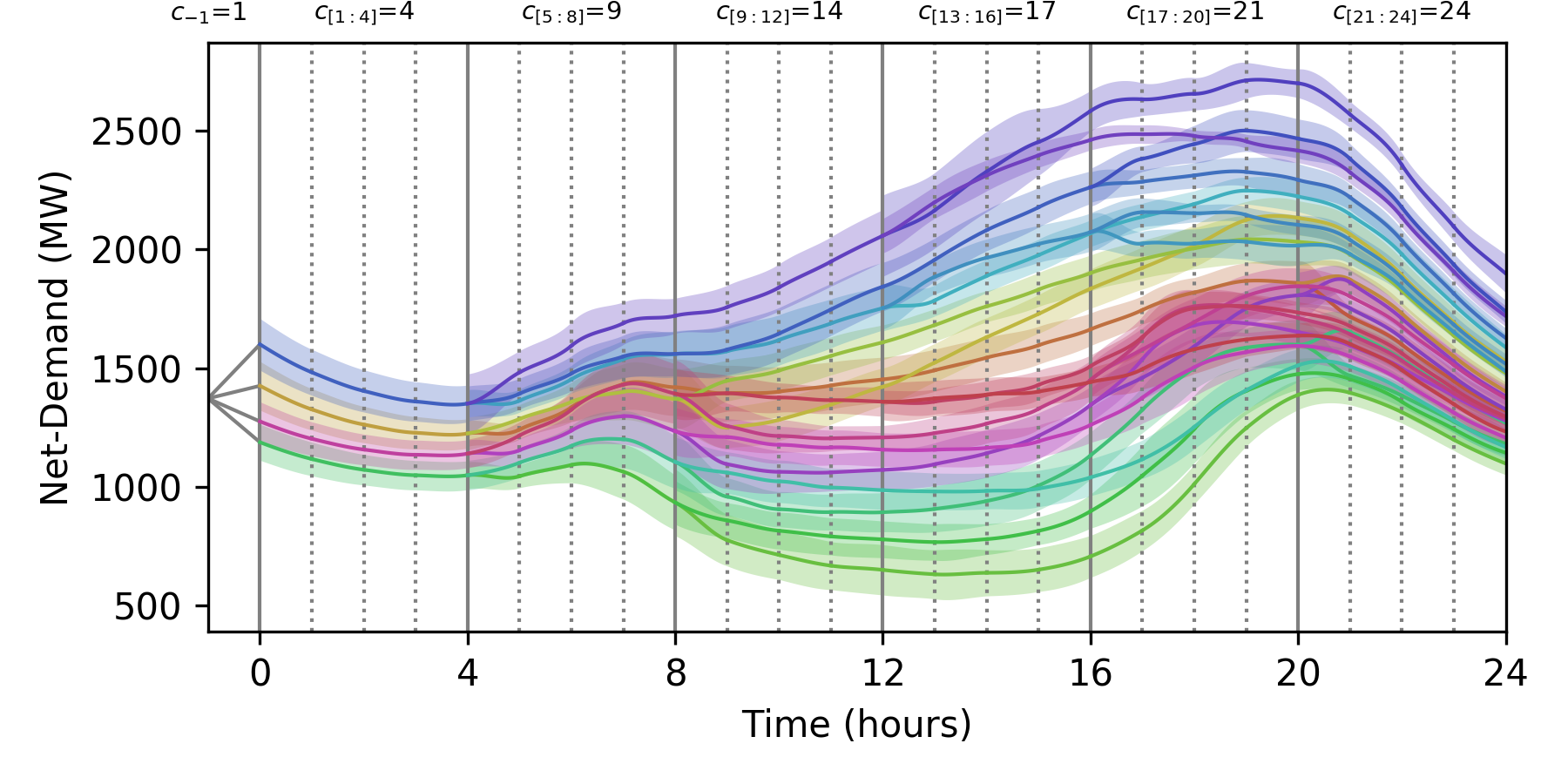}}
  \caption{The zero and third order spline approximation scenario trees (solid colored lines) along with the root mean square error of the scenario bundle yielding the corresponding centroid (shaded region).
The $c_h$ indicator on top shows how many nodes are defined for each stage.}\label{fig:trees}
  \vspace*{-1em}
\end{figure}

To illustrate the proposed framework we use $L=1122$ sample days of \emph{scaled} and \emph{aggregated} CAISO net-load measurements spanning the years 2013-2017, with a 5-minutes resolution\cite{oasis}.
We find the spline approximations of the load scenarios, where each polynomial spans a single hour of load. 
We make two approximations: (a) cubic spline with $C^1$ continuity forced across the nodes and for comparison (b) zero order (discontinuous piece wise constant) spline. The zero order polynomial reflects standard unit commitment practice, where the load is averaged over each hour.
To solve for the spline coefficients we use a standard least-squares approach. 
As an example, Figure \ref{fig:splines} captures these approximations for three sample trajectories.
We split the input scenarios (days) into two sets, 70\% we use for the tree construction and the remaining 30\% we use as sample days for verifying the solution.
The tree is constructed for a ``generic'' load day using a simple recursive $\kmeans$ algorithm which solves in a matter of seconds, but could equally be generated from a more targeted dataset, such as a set of load predictions based on weather forecasts and historical load data.
For both zero and third order spline approximations the training scenarios are reduced down to the trees shown in Figure \ref{fig:trees}. 
The figure shows the root (arbitrarily chosen to be at hour $-1$), the knots of the polynomials (dotted vertical lines), the allowed branching points (solid vertical lines) and the number of nodes per stage indicated by the $c_h$ value on top. 
The different colored paths show the reduced sample path centroids, with the root mean square error indicated by the surrounding shaded region.

\subsection{Solving the Stochastic Unit Commitment}
For generation data we use a single area of the IEEE RTS-96 test case \cite{grigg1999ieee}, but consider it to be a single-bus system.
As mentioned before, the CAISO load data is scaled by a factor of $16$ in order to fit into the range of generation available in the RTS system.
The cost terms used in the objective function \eqref{eq:objective} are summarized in Table \ref{tbl:rtscost} for the different RTS case generation types. 
These terms are based on the original case data, however, a piece wise linear objective is assumed and the quadratic component of the generator cost curves are ignored. 
Fixed capacity costs ($\overline{R}$, $\underline{R}$, $\overline{Y}$) are proportional to the linear and constant cost terms.
\begin{table}
  \caption{Cost terms of the different RTS case generators.}\label{tbl:rtscost}
  \centering
  \begin{tabular}{r|rrrrrrr}
    \toprule
      Unit Type & $\overline{S}$ & $\underline{S}$ & $Y$ & $X$ & $\overline{R}$ & $\underline{R}$ & $\overline{Y}$\\
    \midrule
      U12 & 1500 & 0 &  86.4 &  56.6 &  28.3 &  28.3 &  21.6\\
      U20 & 1500 & 0 & 400.7 & 130.0 &  65.0 &  65.0 & 100.2\\
      U50 & 1500 & 0 &   0.0 &   0.0 &   0.0 &   0.0 &   0.0\\
      U76 & 1500 & 0 & 212.3 &  16.1 &   8.0 &   8.0 &  53.1\\
     U100 & 1500 & 0 & 781.5 &  43.7 &  21.8 &  21.8 & 195.4\\
     U155 & 1500 & 0 & 382.2 &  12.4 &   6.2 &   6.2 &  95.6\\
     U197 & 1500 & 0 & 832.8 &  48.6 &  24.3 &  24.3 & 208.2\\
     U350 & 1500 & 0 & 665.1 &  11.8 &   5.9 &   5.9 & 166.3\\
     U400 & 1500 & 0 & 395.4 &   4.4 &   2.2 &   2.2 &  98.8\\
    \bottomrule
  \end{tabular}
  \vspace*{-1.5em}
\end{table}

Using a net-load tree as the input, we solve the unit commitment problem as described by \eqref{eq:continuity1}-\eqref{eq:objective}, for both the cubic spline approximation (CT-MSRUC) and the discrete time formulation (DT-MSRUC).
The latter is simply a special case of the described formulation, where with respect to the constraints we can imagine the load and generation to be first-order $C^0$ continuous, thus enforcing the ramping constraints and having the same number of free variables as a conventional zero order formulation.
Modeling was done using Python \cite{numpy,scipy,matplotlib,scikit-learn}, and for solving we used Gurobi 7.5 \cite{gurobi} with a MIP gap parameter of $5\%$. 
The solution time of the DT-MSRUC problem was 40 minutes while for the CT-MSRUC it solved in 67 hours. The increase in computational time is expected because the number of continuous variables is roughly quadrupled, while the number of integer variables remains the same. 
Reducing the solver time drastically through decomposition and parameter tuning should be relatively straightforward, but was not the focus of this work. As an example, using a $\rho = 2.5$ solved in a couple of hours.

\begin{figure}
  \centering
  \subfloat[DT-MSRUC]{\includegraphics[width=0.49\columnwidth]{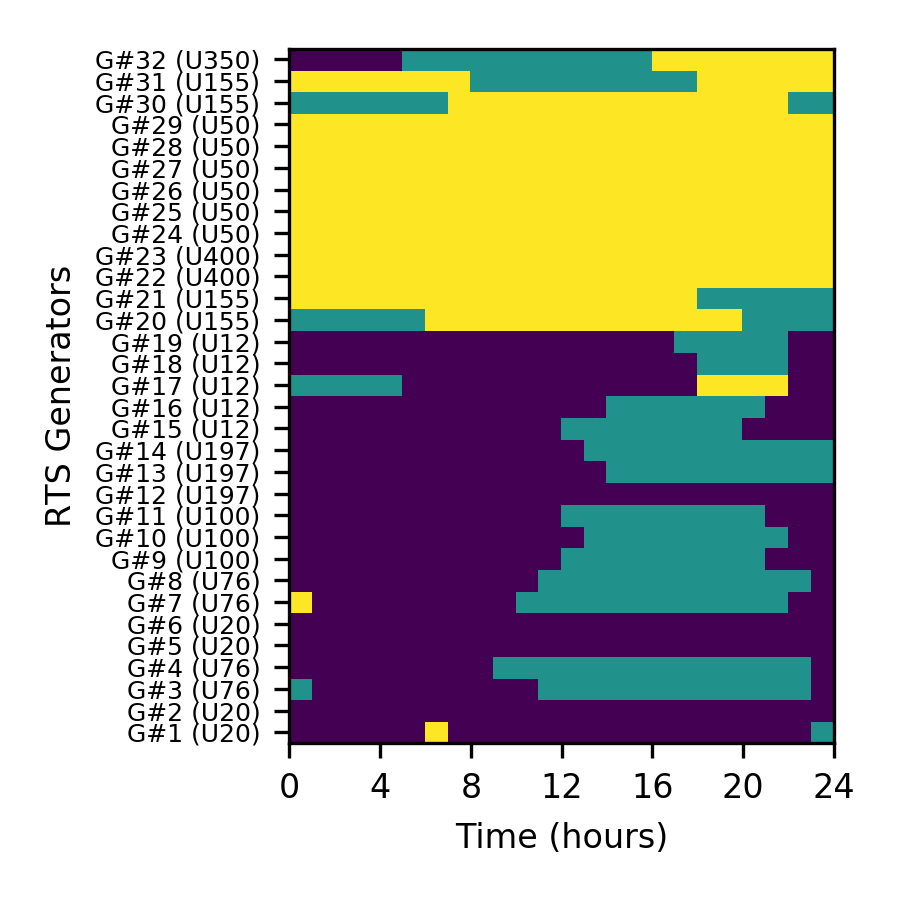}}
  \subfloat[CT-MSRUC]{\includegraphics[width=0.49\columnwidth]{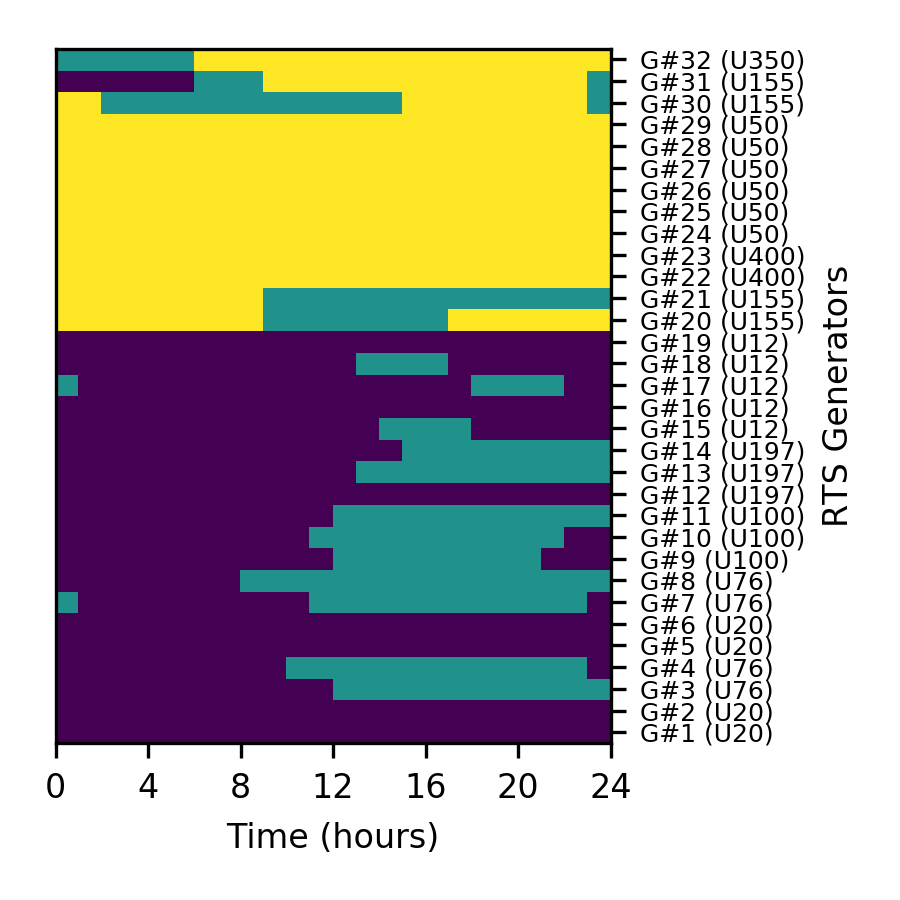}}
  \caption{Unit commitment solutions for the DT-MSRUC and CT-MSRUC formulations. Purple (dark) region denotes units that are never committed, green (gray) denotes units that are committed for one or more paths of the solution, while yellow (light) denotes units committed for the scheduled solution.}\label{fig:commitment}
  \vspace*{-1em}
\end{figure}
We solved the two unit commitment problems CT-MSRUC and DT-MSRUC with the reserve parameter $\rho=3$ and applied the solutions to the set of test scenarios, considering the real-time cost and potential infeasibility issues.
Figure \ref{fig:commitment} visualizes the unit commitment of both solutions. 
The purple shows uncommitted units for any path during a particular hour and yellow shows committed units for the \emph{scheduled} (most likely) path, while green shows what is not committed for the most likely load, but is committed for one or more of the alternative load trajectories. 
We observe that the commitment is fairly similar between the two cases, the hydro (U50) and nuclear (U400) units are committed for most paths and a handful of oil and coal plants (U12, U76, U100, U197) are brought online for the afternoon hours of several paths. A notable difference lies in the usage of the large coal plants (U155, U350) which are scheduled differently for the two solutions.

\begin{figure}
  \centering
  \includegraphics[width=0.95\columnwidth]{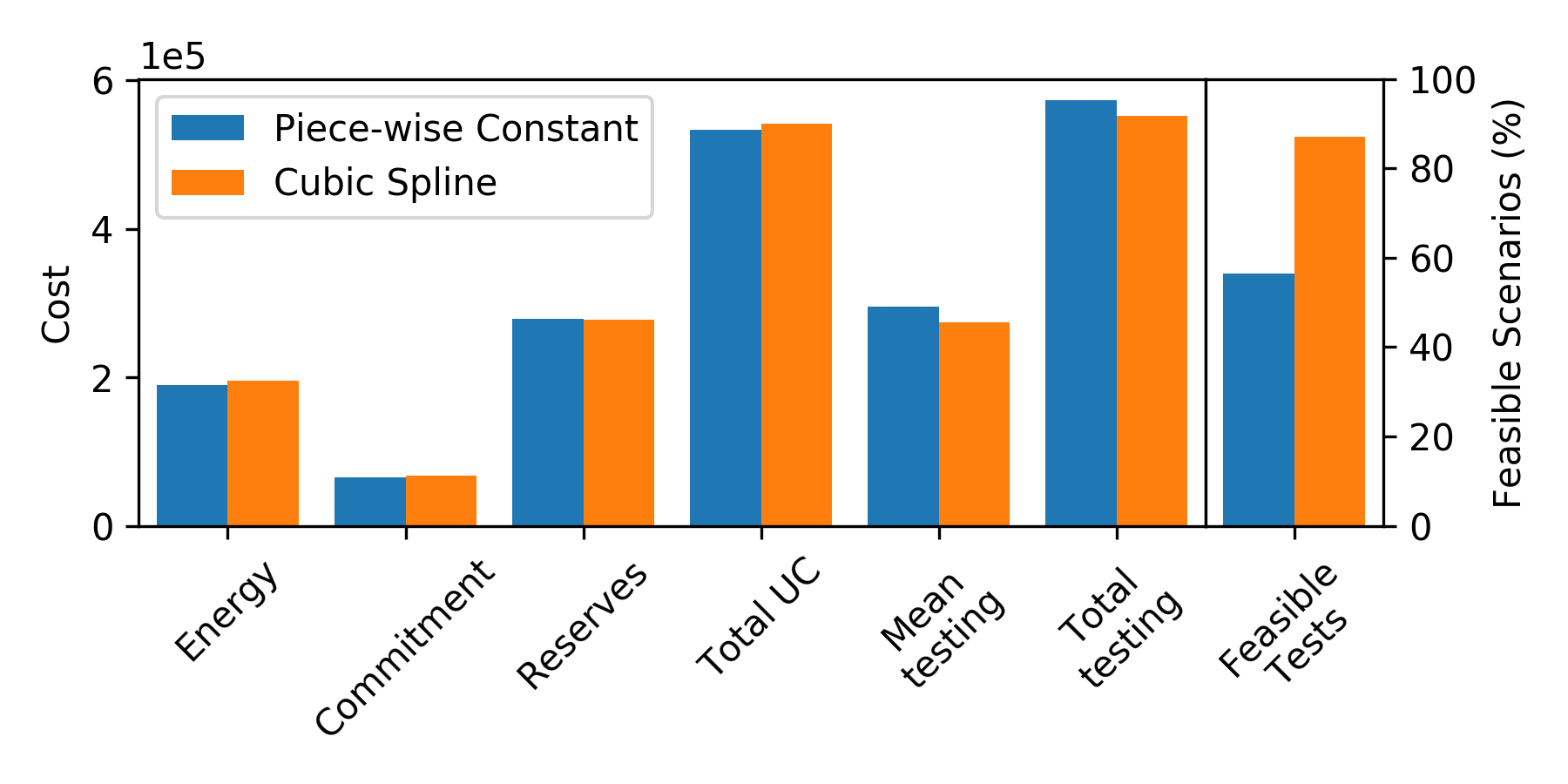}
  \caption{The various cost components of the solution and the test scenarios.}\label{fig:cost}
  \vspace*{-1em}
\end{figure}

Figure \ref{fig:cost} compares several aspects of the two problems. First, the expected energy, commitment and reserve costs as formulated in \eqref{eq:objective} are very similar, with the continuous solution giving a slightly higher expected cost. 
However, when we run through the set of test scenarios the opposite seems true, the mean cost of commitment and energy turns out to be lower (\emph{mean testing} bar) and when we add the reserve cost (\emph{total testing}) we see that the CT-MSRUC solution is actually cheaper, as it offers a more efficient way to service the realistic test trajectory which is continuous in nature.
Moreover, 44\% of the testing load trajectories are simply not feasible for the DT-MSRUC solution, given the commitment and reserve capacity of scheduled units, compared to a level of 14\% infeasibility for the CT-MSRUC solution.

\begin{figure}
  \centering
  \subfloat[Applying the DT-MSRUC solution.]{\includegraphics[width=0.9\columnwidth]{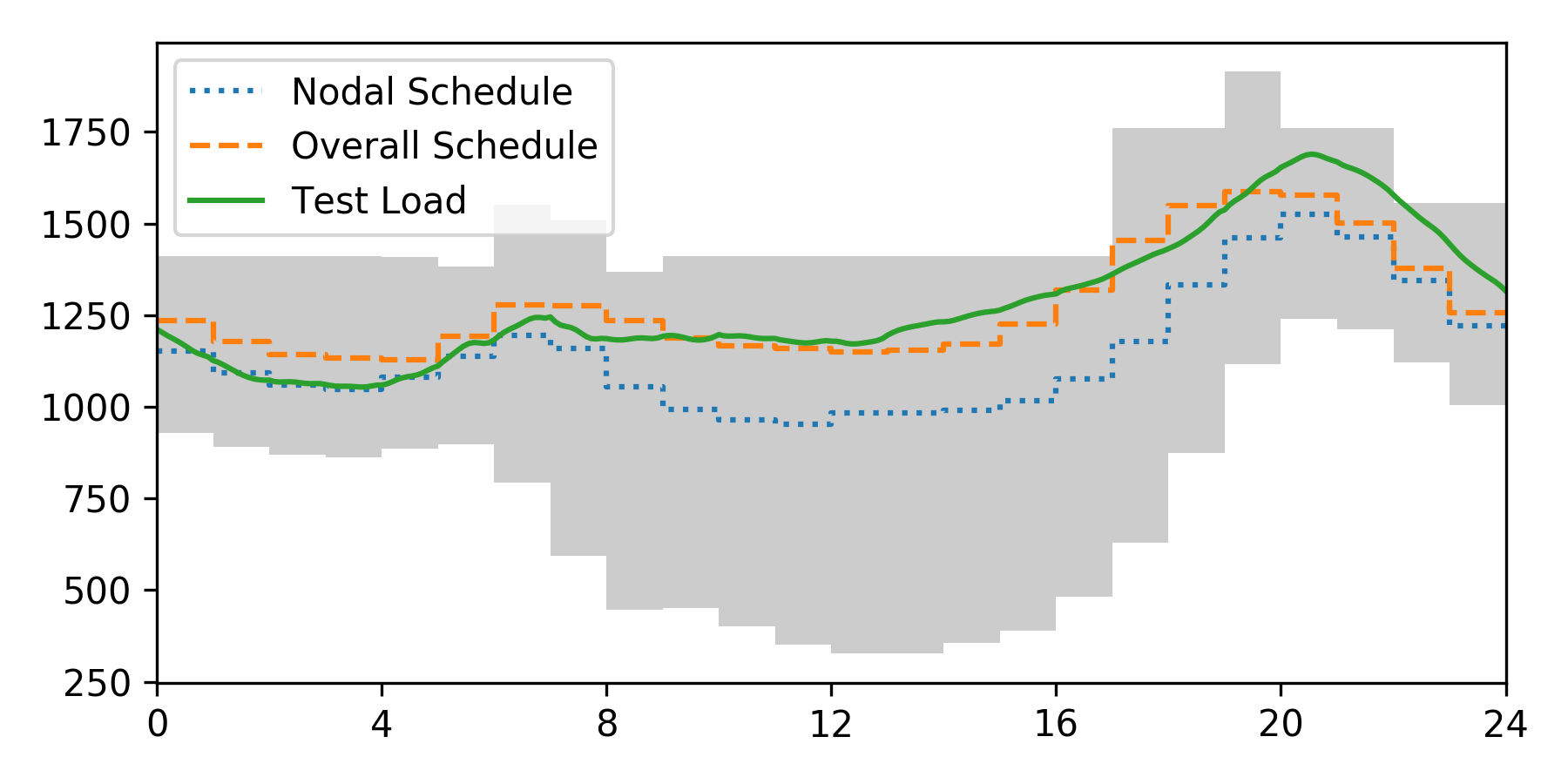}}\\
  \subfloat[Applying the CT-MSRUC solution.]{\includegraphics[width=0.9\columnwidth]{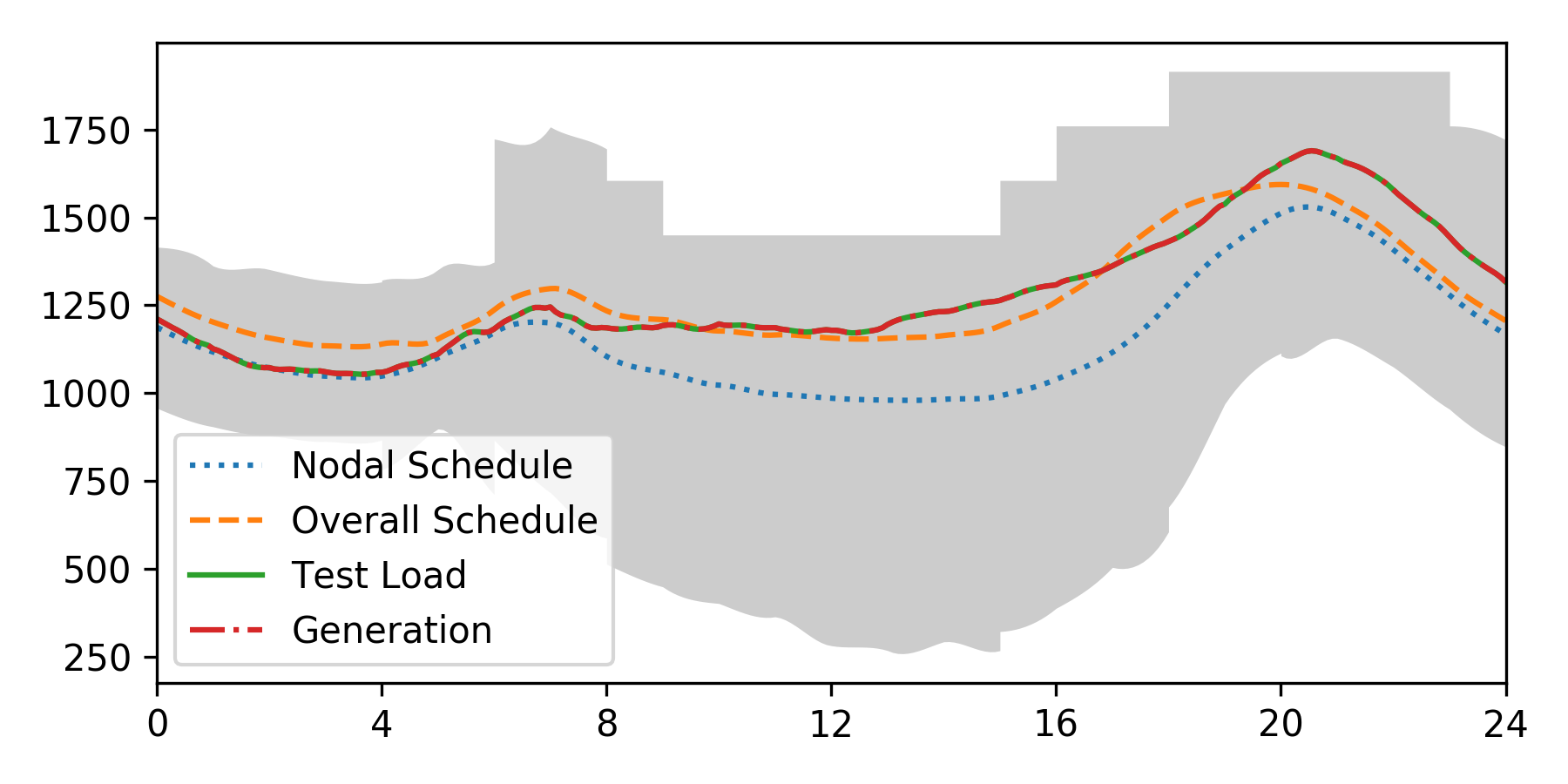}}
  \caption{Generator schedule (dashed orange) along with the dispatch (dotted blue), commitment and reserve solution given along the path of the load scenario tree that closest resembles a given test load scenario (solid green). The reserve capacity of committed units is denoted by the shaded gray region.}\label{fig:tests}
  \vspace*{-1.5em}
\end{figure}

In real-time, the system operator must maintain generation--consumption balance. As the load prediction is never perfect, the exact dispatch solution of the UC is not followed, but the commitment and reserve capacity "schedule" should be known ahead of real-time.
For the case of stochastic UC, however, each path of the scenario tree will contain a different commitment schedule. 
Observing the load in real-time, we follow the tree to obtain the specific commitment profile, essentially throwing out any diverging paths until we end up with what looks like a solution to a deterministic UC.
Figure \ref{fig:tests} shows a given test load trajectory (solid green) and the nearest load path on the scenario tree (dotted blue), compared to the most likely path on the tree (dashed yellow).
While the load will not exactly match the branch of the tree we follow, we do know the commitment, scheduled dispatch and reserve capacity, and can transform these into the range of load that our solution can service (shaded gray region of figure).
Here we can observe the CT-MSRUC solution having two major advantages over the DT-MSRUC solution, leading to higher likelihood of the solution being feasible for a real load trajectory. First, we do a much better job of incorporating realistic ramps, so we can better anticipate load trajectories that change abruptly and have the right fleet of units available, and secondly, by not modeling in terms of averaged-load of an hour, we are less likely to "cut corners"  where the test load trajectory falls outside the feasible region for a brief time. This can be seen in Figure \ref{fig:tests} where the discrete-time solution can not deliver desired load during the first minutes of hour 22.

\section{Conclusions \& Future Work}
This paper describes a stochastic continuous-time unit commitment formulation, based on Bernstein polynomials.
We briefly discuss how a set of input scenarios are reduced into a scenario tree, and formulate a multi-stage stochastic unit commitment, where commitment of units varies between the different paths of the tree.
We show that our formulation outperforms conventional stochastic discrete-time unit commitment, as we better account for necessary real-time ramps and we are less likely to underestimate reserves needed for the hours where load is sharply increasing or decreasing.

Future work will go more in depth into the tree construction, focus on extending this to multi-bus systems, incorporating the spatial correlation of load into the scenario tree, as well as looking at the trade-off between tree sizes, branching structure and higher order polynomials. 
Further, a natural progression of this work is incorporating decomposition methods into the solution algorithm to compensate for the unavoidable increase in complexity.

\bibliography{MainReferences}
\bibliographystyle{IEEEtran}
\end{document}